\documentclass[aos,preprint,12pt]{imsart}
\RequirePackage[OT1]{fontenc}

\RequirePackage{graphicx,amsthm,amsmath,natbib,lineno} 
\RequirePackage[colorlinks,citecolor=blue,urlcolor=blue]{hyperref} 
\RequirePackage{hypernat} 
\usepackage{bbm}
\startlocaldefs 
\newtheorem{definition}{Definition}[section] 
\newtheorem{thm}{Theorem}[section] 
 
\newtheorem{corollary}{Corollary}[section] 
\newtheorem{example}{Example}[section] 
 
\newtheorem{Algorithm}{Algorithm}[section] 
\newtheorem{Remark}{Remark}[section] 
\newtheorem{Notation}{Notation}[section] 
\endlocaldefs 

\begin{document} 

\begin{frontmatter}
  \title{Minimal Markov Models
    \protect\thanksref{T1}} \runtitle{Minimal Markov Models } \thankstext{T1}{This work is partially supported by
    PRONEX/FAPESP Project Stochastic behavior, critical phenomena and
    rhythmic pattern identification in natural languages (grant number
    03/09930-9) and by CNPq Edital Universal (2007), project:
    ``Padr\~{o}es r\'{ \i}tmicos, dom\'{ \i}nios pros\'{o}dicos e
    modelagem probabil\'{ \i}stica em corpora do portugu\^{e}s''.}
\begin{aug}
	\author{\fnms{Jes\'{u}s} \snm{ E. Garc\'{ \i}a}\thanksref{t1,m1}\ead[label=e1]{jg@ime.unicamp.br}} \and 
	\author{\fnms{Ver\'{o}nica} \snm{ A. Gonz\'{a}lez-L\'{o}pez}\thanksref{t1,m1}\ead[label=e2]{veronica@ime.unicamp.br}}
	
	\thankstext{t1}{Departamento de Estat\'{ \i}stica. Intituto de
          Matem\'{a}tica Estat\'{ \i}stica e Computa\c c\~{a}o
          Cient\'{ \i}fica.}
	
	\runauthor{J. E. Garc\'{ \i}a and V. A. Gonz\'{a}lez-L\'{o}pez}
	
	\affiliation{Universidade Estadual de Campinas \thanksmark{m1}}
	
	\address{Departamento de Estat\'{ \i}stica\\
	Instituto de Matem\'{a}tica Estat\'{ \i}stica e Computa\c c\~{a}o Cient\'{ \i}fica\\
	Universidade Estadual de Campinas \\
	Rua Sergio Buarque de Holanda,651\\
	Cidade Universit\'{a}ria-Bar\~{a}o Geraldo\\
	Caixa Postal: 6065\\
	13083-859 Campinas, SP, Brazil\\
	\printead{e1}\\
	\phantom{E-mail:\ }\printead*{e2}}
	
\end{aug}
\begin{abstract}
	In this work we introduce a new and richer class of finite order Markov chain models and address the following model selection problem: find the Markov model with the minimal set of parameters (minimal Markov model)  which is necessary to represent a source as a Markov chain of finite order. Let us call $M$ the order of the chain and $A$ the finite alphabet, to determine the minimal Markov model, we define an equivalence relation on the state space $A^{M}$, such that all the sequences of size $M$ with the same transition probabilities are put in the same category. In this way we have one set of $(|A|-1)$ transition probabilities for each category, obtaining a model with a minimal number of parameters. We show that the model can be selected consistently using the Bayesian information criterion. 
\end{abstract}
\end{frontmatter}

\section{Introduction}

 In this work we consider discrete stationary processes over a finite alphabet $A.$ Markov chains of finite order are widely used to model stationary processes with finite memory. A problem  with full Markov chains models of finite order $M$ is that the number of parameters $(|A|^M(|A|-1))$ grows exponentially with the order $M,$ where $|A|$ denotes the cardinal of the alphabet $A.$ Another characteristic is that the class of full Markov chains is not very rich, fixed the alphabet $A$ there is just one model for each order $M$ and in practical situations could be necessary a more flexible structure in terms of number of parameters. For an extensive  discussion of those two problems se \citet{Buhlmann1999}. A richer class of finite order Markov models introduced by \citet{Rissanen1983} and \citet{Buhlmann1999} are the variable length Markov chain models (VLMC) which are mentioned in section \ref{GPCT}. In the VLMC class, each model is identified by a prefix tree $\cal T$ called context tree. For a given model with a context tree $\cal T$, the final number of parameters for the model is $|{\cal T}|(| A|-1)$ and depending on the tree, this produce a parsimonious model. In \citet{Csiszar2006} is proved that the bayesian information criterion (BIC) can be used to consistently choose the VLMC model in an efficient way using the context tree weighting (CTW) algorithm.

In this paper we introduce a larger class of finite order Markov models, and we address the problem of model selection inside this class, showing that the model can be selected consistently using the BIC criterion. In our class, each model is determined by choosing a partition of the state space, our class of models include the full Markov chain models and the VLMC models because a context tree can be seen as a particular partition of the state space (see for illustration the example \ref{CTpartition}).

In  Section \ref{minimalmodels},  we define the minimal Markov models and show that this models can be selected in a consistently in theorems  \ref{GeralREFTEO1} and  \ref{GeralREFTEO2}.
In Section \ref{algoritmos} we show two algorithms that use the results in Section \ref{minimalmodels} to choose consistently a minimal Markov model for a sample and some simulations. Section \ref{conclu} have the conclusions and Section \ref{demostra} have the proofs.

 \section{Minimal Markov models}\label{minimalmodels}

\subsection{Notation}  ${  }$
Let  $(X_t)$ be a discrete time order $M$ Markov chain on a finite alphabet  $A$.  Let us call ${\cal S}=A^M$ the state space. Denote the string $a_ma_{m+1} \ldots a_{n}$ by $a_m^n,$ where  $\,a_i \in A,\,m\leq i \leq n.$

Let ${\cal L}=\{ L_1,L_2,\ldots, L_K\}$ be a partition of ${\cal S},$
 \begin{eqnarray}
 P(L,a)=\sum_{ s \in L}\mbox{Prob}(X_{t-M}^{t-1}=s,X_t=a), \,\,\, a \in A, \,\, L \in {\cal L};
 \end{eqnarray}
 \begin{eqnarray}
 P(L)= \sum_{ s \in L} \mbox{Prob}(X_{t-M}^{t-1}=s),\,\,\, L \in {\cal L}.
 \end{eqnarray}
Let $x_1^n$ be a sample of the process $\big( X_t \big),\,s \in {\cal S},$  $a \in A$ and $n>M.$ We denote by $N_n(s,a)$ the number of occurrences of  the string $s$  followed by $a$ in the sample $x_1^n,$
\begin{eqnarray} \label{contaNsa}
N_n(s,a)= \big|\{t:M< t \leq n, x_{t-M}^{t-1}=s ,x_{t}=a \}\big|, 
\end{eqnarray} 
the number of occurrences of $s$ in the sample $x_1^n$ is denoted by $N_n(s)$ and
\begin{eqnarray}
N_n(s)=\big|\{t:M< t \leq n, x_{t-M}^{t-1}=s \}\big|.
\end{eqnarray}

The number of occurrences of elements into $L$ followed by $a$ is given by,
\begin{eqnarray}\label{defi1}
N^{\cal L}_n(L,a)=\sum_{s \in L} N_n(s,a),\,\, L \in {\cal L};
\end{eqnarray}
the accumulated number of $N_n(s)$ for $s$ in $L$ is denoted by, 
\begin{eqnarray}\label{defi2}
N^{\cal L}_n(L)=\sum_{s \in L} N_n(s),\,\, L \in {\cal L}.
\end{eqnarray}

\subsection{Good partitions of $\cal S$}\label{pbuenas}  ${  }$

\begin{definition}
	Let  $(X_t)$ be a discrete time order $M$ Markov chain on a finite alphabet  $A,$ ${\cal S}=A^M$ the state space. A partition ${\cal L}=\{ L_1,L_2,\ldots, L_K\}$ of $\cal S$ is a good partition of $\cal S$ if for each $s, s^{\prime}\,\in L, \,\,\,L \in {\cal L},$ \[Prob(X_{t}=. \,|X^{t-1}_{t-M}=s)=Prob(X_{t}=.\,|X^{t-1}_{t-M}=s^{\prime}).\]
\end{definition}

\begin{Remark} For a discrete time order $M$ Markov chain on a finite alphabet  $A$ with ${\cal S}=A^M$ the state space, ${\cal L}={\cal S}$  is a good partition of ${\cal S}.$
\end{Remark}

If  ${\cal L}$ is a good partition of ${\cal S},$ we define for each category $L \in {\cal L}$
\begin{eqnarray}\label{goodpartitions}
P(a \vert L)=\mbox{Prob}(X_{t}=a|X_{t-M}^{t-1}=s) \,\,\, \forall a\in A,
\end{eqnarray}
where $s$ is some element into $L.$
As a consequence, if we write $P(x_1^n)=\mbox{Prob}(X_1^n=x_1^n),$ we obtain
\begin{eqnarray}\label{probsample2}
 P(x_1^n)=P(x_1^M)\prod_{L \in {\cal L}, a \in A} P(a \vert L)^{N^{\cal L}_n(L,a)}.
 \end{eqnarray}
 In the same way that \citet{Csiszar2006} we will define our BIC criterion using a modified maximum likelihood.
 We will call maximum likelihood to the maximization of the second term in the equation (\ref{probsample2}) for the given observation. For the sequence  $x_1^n,$ will be 
 \begin{eqnarray}\label{MLgoodp}
 \mbox{ML}({\cal L},{x_1^n})=\prod_{L \in {\cal L}, a \in A} \left( \frac{r_n(L,a)}{r_n(L)} \right)^{N_n^{\cal L}(L,a)},
 \end{eqnarray}
 where 
 
 \begin{eqnarray}\label{erres}
r_n\big(L,a\big)=\frac{N^{{\cal L}}_n(L,a)}{n}, \,\,\,\,\,\, a \in A, \,\,\, L \in {\cal L}\,\,&and&
r_n\big(L\big)=\frac{N^{{\cal L}}_n(L)}{n}, \,\,\,\, L \in {\cal L}.
\end{eqnarray}
 
The  BIC is given by the next definition
\begin{definition}\label{BICML}
Given a sample $x_1^n,$ of the process $(X_t),$ a discrete time order $M$ Markov chain on a finite alphabet  $A$ with ${\cal S}=A^M$ the state space and ${\cal L}$ a good partition of ${\cal S}.$ The BIC of the model (\ref{MLgoodp}) is given by
\begin{eqnarray}
BIC({\cal L}, x_1^n)=\ln\left(\mbox{ML}({\cal L},x_1^n)\right)-\frac{(|{A}|-1)|{\cal L}|}{2}\ln(n). \nonumber
\end{eqnarray}
\end{definition}

\subsection{Good partitions and context trees} \label{GPCT}  ${  }$

Let  $(X_t)$ be a finite order Markov chain taking values on $A$ and {\sl ${\cal T}$} a set of sequences of symbols from $A$ such that no string in ${\cal T}$ is a suffix of another string in ${\cal T},$ for each $s\in \cal T,$ $d({\cal T})=\max\big(l(s), s \in {\cal T}\big)$ where $l(s)$ denote the length of the string $s,$ with $l(\emptyset)=0$ if the string is the empty string.
\begin{definition}
	${\cal T}$ is a context tree for the process  $(X_t)$ if for any sequence of symbols in $A$, $x^n_1$ sample of the process with $n\geq d({\cal T}),$ there exist $s \in \cal T$ such that$$Prob(X_{n+1}=a|X_1^n=x_1^n)=Prob(X_{n+1}=a|X_{n-l(s)+1}^n=s)$$
\end{definition}
$d({\cal T})$ is the  depth of the tree. \\
The context tree is the minimal state space of the variable length Markov chain (VLMC), \citet{Buhlmann1999}. The context tree for a VLMC with finite depth $M$ define a good partition on the space ${\cal S}=A^M$ as illustrated by the next example.
\begin{example} \label{CTpartition} Let be a VLMC
 over the alphabet $A=\{0,1\}$ with depth $M=3$ and contexts,
$$\{ 0\}, \{01\}, \{ 011\}, \{ 111\} $$
This context tree correspond to the good partition
$\{ L_1, L_2, L_3, L_4 \}$ where\\
  $L_1=\{\{ 000\},\{ 100\},\{ 010\}, \{110\}\},\,L_2=\{\{001 \}, \{101\}\},\,\,L_3=\{011\}$ and $L_4=\{ 111\}. $
\end{example}

\subsection{Smaller good partitions}

\begin{definition}\label{juntarLiLj}
Let ${\cal L}^{ij}$ denote the partition 
\begin{eqnarray}
 {\cal L}^{ij} = \{L_1, \ldots,L_{i-1},L_{ij},L_{i+1},\ldots,L_{j-1},L_{j+1},\ldots, L_K\},\nonumber
 \end{eqnarray}
where ${\cal L} = \{L_1, \ldots, L_K\}$ is a good partition of $\cal S,$ and for $1\leq i<j\leq K$ with $ L_{ij}=L_i\cup L_j.$
\end{definition}
Now we adapt the notation established for the partition ${\cal L}$ to the new partition ${\cal L}^{ij} .$\\

\begin{Notation}\label{Notationconteoij}
for $a \in A$ we  write, 
 \begin{eqnarray}
 P(L_{ij},a)&=&P(L_i,a)+P(L_j,a)\nonumber;\\
 P(L_{ij})&=& P(L_i)+P(L_j) \nonumber.
 \end{eqnarray}
  
\begin{eqnarray}\label{nlij} 
N_n^{{\cal L}^{ij}}(L_{ij},a)= N^{\cal L}_n(L_i,a)+N^{\cal L}_n(L_j,a);
\end{eqnarray}
\begin{eqnarray}
N_n^{{\cal L}^{ij}}(L_{ij})= N^{\cal L}_n(L_i)+N^{\cal L}_n(L_j);
\end{eqnarray}
\end{Notation} 

If $P(. \vert L_i)=P(. \vert L_j)$ then ${\cal L}^{ij}$ is a good partition and (\ref{goodpartitions}) remains valid for ${\cal L}^{ij},$ just is necessary to change ${\cal L}$ by ${\cal L}^{ij}$ in equations (\ref{probsample2}), (\ref{MLgoodp}) and definition (\ref{BICML}). \\
In the following theorem, we show that the BIC criterion provides a consistent way of detecting smaller good partition.

\begin{thm}  \label{GeralREFTEO1}  Let $(X_t)$ be a Markov chain with order $M$ over a finite alphabet $A,$ ${\cal S}=A^{M}$ the state space.
If ${\cal L}=\{L_1, L_2, \ldots , L_K \}$ is a good partition of ${\cal S}$ and ${ L}_i \neq { L}_j,\, L_i, L_j \in {\cal L}.$ Then, eventually almost surely as $n \to \infty,$
 $$  I_{\{BIC({\cal L}^{ij},x_1^n)>BIC({\cal L },x_1^n) \}}=1$$ if, and only if $$P(a|{ L}_i)=P(a|{ L}_j) \; \forall a \in A.$$ 
Where $I_{A}$ is the indicator function  of $A,$ and the ${\cal L}^{ij}$ partition is defined under ${\cal L}$ by equation (\ref{juntarLiLj}). 
\end{thm}
 
 Next we extract from the previous theorem the relation that we use in the next section, in practice to find smaller good partitions. 
 \begin{definition} Let be $(X_t)$ a Markov chain of order $M,$ with finite alphabet $A$ and state space   ${\cal S}=A^M$, $x_1^n$ a sample of the process and let  ${\cal L}=\{ L_1,L_2,\ldots, L_K\}$ be a good partition of ${\cal S},$
  \begin{eqnarray} 
  d_{\cal L}(i,j)= \frac{1}{\ln(n)} \sum_{ a \in A}\left\{
N^{{\cal L}}_n(L_i,a)\ln\left(\frac{N_n(L_i,a)}{N_n(L_i)}\right)\right. + N^{{\cal L}}_n(L_j,a)\ln\left(\frac{N_n(L_j,a)}{N_n(L_j)}\right) \nonumber \\
\left.  -N^{{\cal L}^{ij}}_n(L_{ij},a)\ln\left(\frac{N_n(L_{ij},a)}{N_n(L_{ij})}\right)\right\} 
\end{eqnarray}
   
 \end{definition}

\begin{corollary}\label{GeralREFEQ1} Under the assumptions of theorem \ref{GeralREFTEO1},

 \begin{eqnarray}
BIC({\cal L}, x_1^n)-BIC({{\cal L}^{ij}}, x_1^n)<0 &\iff&  d_{\cal L}(i,j)<\frac{(|{A}|-1)}{2}. \nonumber
\end{eqnarray}
\end{corollary}
\begin{proof} From equation (\ref{dgeral}) we have the validity of the result.
\end{proof}
 
 \begin{Remark}The results will remain valid if we replace the constant $\frac{(|A|-1)}{2}$ for some arbitrary constant, positive and finite value $v,$ into the definition (\ref{BICML}). 
\end{Remark}

 \begin{Remark}  \label{GOODBADTEO} Under the assumptions of theorem \ref{GeralREFTEO1},
if  $P(a \vert L_i) \neq P(a \vert L_j)$ for some $a \in A,$  then 
eventually almost surely as $n \to \infty,$
 $  BIC({\cal L},x_1^n)>BIC({\cal L}^{ij},x_1^n)$
where ${\cal L}^{ij}$ verified the definition (\ref{juntarLiLj}).
\end{Remark}

\subsection{ Minimal good partition}   ${  }$

We want to find the smaller good partition into the universe of all possible good partitions of ${\cal S}.$ This special good partition could be defined as follows and it allows the definition of the most parsimonious model into the class considered in this paper.

\begin{definition}
	Let  $(X_t)$ be a discrete time order $M$ Markov chain on a finite alphabet  $A,$ ${\cal S}=A^M$ the state space. A partition ${\cal L}=\{ L_1,L_2,\ldots, L_K\}$ of $\cal S$ is the minimal good partition of $\cal S$ if, $\forall  L \in {\cal L},$
	 \[s, s^{\prime} \in L\,\,\mbox{  if, and only if}\,\, Prob(X_{t}= .\, |X^{t-1}_{t-M}=s)=Prob(X_{t}= . \,|X^{t-1}_{t-M}=s^{\prime}).\]
\end{definition}

\begin{Remark} For a discrete time order $M$ Markov chain on a finite alphabet  $A$ with ${\cal S}=A^M$ the state space, $\exists !$ minimal good partition of $\cal S.$ 
\end{Remark}
In the next example we emphasize the difference between good partitions and the minimal good partition,

 The next theorem shows that for $n$ large enough we achive the partition ${\cal L}^{*}$ which is the minimal good partition.

\begin{thm}  \label{GeralREFTEO2}   Let $(X_t)$ be a Markov chain with order $M$ over a finite alphabet $A,$ ${\cal S}=A^{M}$ the state space  and let $\cal P$ be the set of all the partitions of $\cal S.$ Define,
 $$ {{\cal L}_{n}^{*}}=argmax_{\cal L \in P} \{BIC({\cal L},x_1^n)\}$$ 
 then, eventually almost surely as $n \to \infty,$
 $$ {\cal L}^{*}={{\cal L}_{n}^{*}}$$
\end{thm}

\section{Minimal good partition estimation algorithm}\label{algoritmos}

\begin{Algorithm}\label{al2}
	({\it MMM} algorithm for good partitions)
	\\
	Consider $x_1^n$ a sample of the Markov process $(X_t),$  with order $M$ over a finite alphabet $A,$ ${\cal S}=A^{M}$ the state space. \\ Let be  $ {\cal L}=\{L_1,L_2,\ldots ,L_K\}$ a good partition of ${\cal S},$ for each $s \in {\cal S},$\\
	\begin{itemize}
	\item[1] for $i=1,2,\cdots,K-1,$
	\begin{itemize} \item[] for $j=i+1,2,\cdots,K,$ 
		\begin{itemize} \item[] Calculate $d_{\cal L}(i,j)$
		\item[]  $R_n^{i,j}= I_{\{d_{\cal L}(i,j)<\frac{(|{A}|-1)}{2}\}}$
		\end{itemize}
		\end{itemize}	
	\item[2]If $R_n^{i,j}=1,$ define $L_{ij}=L_i\cup L_j$ and ${\cal L}={\cal L}^{ij}$ . Else $i=i+1,$ Return to step 1
	\end{itemize}
\end{Algorithm}
The algorithm allows to define the next relation based on the sample $x_1^n,$ 
\begin{definition}
 for  $r, s \in  {\cal S}; \;\; r\sim_n s \iff R_n^{a(r),a(s)}=1.$
\end{definition}

For $n$ large enough, the algorithm return the minimal good partition.
\begin{corollary} Let {\sl $\{X_t, t=0,1,2,\ldots \}$} be a Markov chain with order $M$ over a finite alphabet $A$, ${\cal S}=A^M$ and $x_1^n$ a sample of the Markov process.
	$\hat{\cal L}_n, $ given by the algorithm (\ref{al2}) converges almost surely eventually to ${\cal L}^*,$ where  ${\cal L}^*$ is the minimal good partition of ${\cal S}.$ 
\end{corollary}
\begin{proof}
Because $K <\infty,$ for $n$ large enough, the algorithm return the minimal good partition.
\end{proof}

\begin{Remark}
In the worst case, which correspond to an initial good partition equal to $\cal S$, we need to calculate the term $\left( \frac{N^{{\cal L}}_n(L,a)}{N^{{\cal L}}_n(L)} \right)^{N_n^{\cal L}(L,a)}$ for each $s \in {\cal S}$ plus $K(K-1)/2$ divisions to implement the algorithm (\ref{al2}). 
\end{Remark}

The next algorithm is a variation of the first. In this case the partitions are grow selecting the pair of elements with the minimal value of $\{d_{\cal L}(i,j)$, the algorithm stop when there is not $\{d_{\cal L}(i,j)$ lower than $(|A|-1)/2$.

\begin{Algorithm}\label{almin}
	Consider $x_1^n$ a sample of the Markov process $(X_t),$  with order $M$ over a finite alphabet $A,$ ${\cal S}=A^{M}$ the state space. \\ Let  $ {\cal L}=\{L_1,L_2,\ldots ,L_K\}$ be a good partition of ${\cal S}$
	\begin{itemize}
	\item[1] Calculate  $$ (i^{*},j^{*})=arg \min_{i,j | 1\leq i <j\leq K} \{ d_{ {\cal L}(i,j)} \} $$	
	\item[2] If $ d_{ {\cal L}(i,j)}<\frac{|A|-1}{2}$ then $ {\cal L}= {\cal L}{i^{*}j^{*}},\;\; K=K-1$ and return to 1.\\ Else end.
	\end{itemize}
\end{Algorithm}

This algorithm is consistent and always return a partition but have a greater computational cost. Taking in consideration that the cost depend on $K$ and  that for a Markov chain of order $M$ we consider samples of size $n$ such that $\log(n)>M$.  The two algorithms \ref{al2} and \ref{almin} have a computational cost that is linear in $n$ (the sample size).

\subsection{Dendrograms and  {\it MMM} algorithm } \label{dendrogamas}

In practice, when the sample size is not large enough and the algorithm \ref{al2} has not converged, it is possible that the algorithm will not return a partition of $\cal S,$ independent of the value used in $v.$ In that case, a better approach can be to use for each $r, s \in {\cal S}$  the function   $d_n(r,s)$ as a similarity measure between $r$ and $s$. 
Then $d_n(r,s)$ can be used to produce a dendrogram and then use the partition defined by the dendrogram as the  partition estimator. 

Also in practice it is possible that the maximum number of free parameters in our model is limited by a number $K$. In that case, the logic choice will be to find a value of $d$ in the dendrogram such that the size of the partition obtained cutting the dendrogram in $d$ is less or equal to $K$, the chosen model will be the one defined by that partition.

\begin{example}\label{exmodel}

Consider a Markov chain of order $M=3$ on the alphabet $A=\{0,1,2\}$ with classes:

\begin{eqnarray*}
L_1&=&\{000,100,200,010,110,210,020,120,220,022,122,222 \},\\
L_2&=&\{001,101,201,011,111,211,021,121,221 \},\\
L_3&=&\{012, 112, 212, 002 \},\\
L_4&=&\{102\},\\
L_5&=&\{202\},
\end{eqnarray*}  
and transition probabilities,
\begin{eqnarray*}
P(0|L_1)=0.2, \,\,\,\,P(1|L_1)=0.3,\\
P(0|L_2)=0.4, \,\,\,\, P(1|L_2)=0.3,\\
P(0|L_3)=0.4, \,\,\,\, P(1|L_3)=0.1,\\
P(0|L_4)=0.1, \,\,\,\, P(1|L_4)=0.4,\\
P(0|L_5)=0.3, \,\,\,\, P(1|L_5)=0.5.
\end{eqnarray*}  

On this example, $|A|=3$ so the penalty constant is $1=\frac{|A|-1}{2}$.
We simulated samples of sizes  $n=5000$ and $9000$, obtaining dendrograms on figure \ref{dendrogramas}. The dendrogram for the sample size of $9000$ gives the correct partition.

\begin{figure}
	[h!] 
	
	\label{dendrogramas}
		
		\centering
	\includegraphics[width=3in]{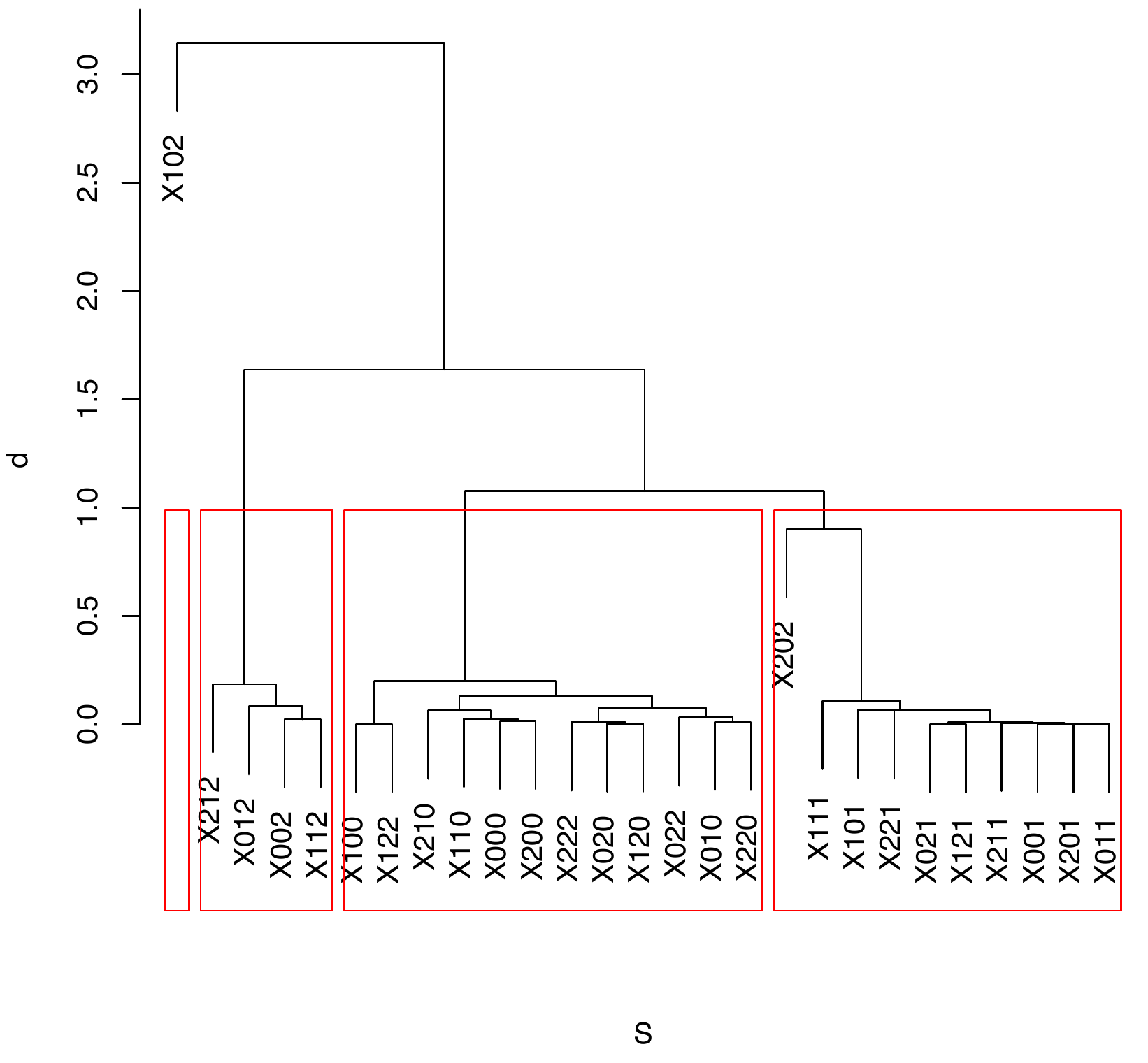}
			\includegraphics[width=3in]{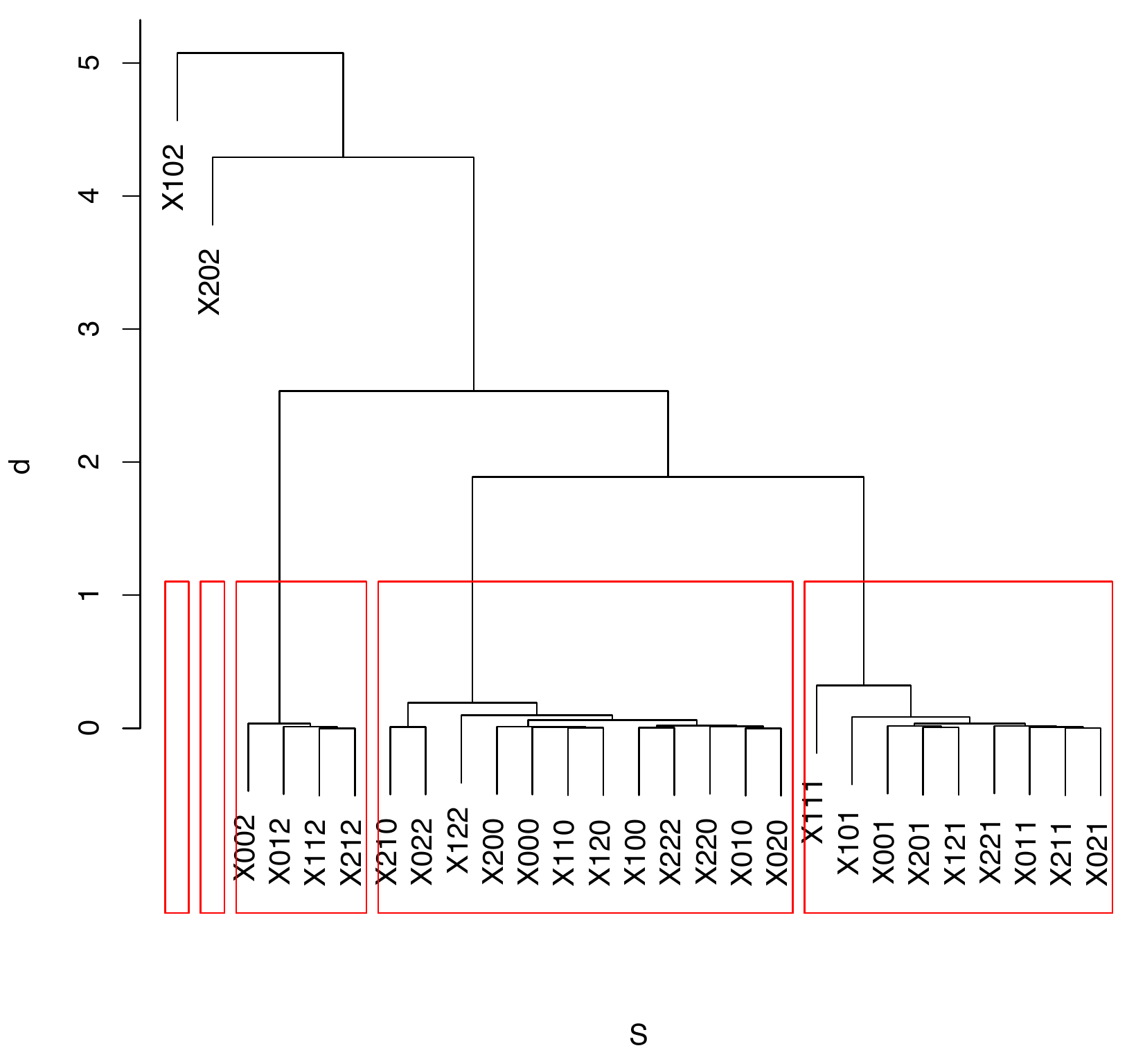}
			
			\caption{The figure shows the dendrograms for the model on example \ref{exmodel}  estimated using algorithm  \ref{REFALG1} for sample sizes of $5000$ (upper picture) and $9000$ (lower picture).}
	\end{figure}

\end{example}

\subsection{Simulations}\label{simulacion1}

We implemented a simulation study for the model described on example  \ref{exmodel}. More precisely we simulated $1000$ samples of the process for each of the sample sizes $4000, 6000, 8000$ and $10000$. For each sample we calculate the values $d_n(r,s)$ and build the corresponding dendrogram (using the R-project package hclust with linkage method complete). Table \ref{tablaerr1} show the results.

\begin{table}
[h!]
\caption{Number of errors on the partition estimated for the model on example \ref{exmodel}}\label{tablaerr1}
\begin{center}
\begin{tabular}{c|c}
Sample size & Proportion of errors\\
\hline
4000 & 0.801\\
6000 & 0.495\\
8000 & 0.252\\
10000 & 0.161\\
\end{tabular}
\end{center}
\label{default}
\end{table}%

\subsection{Simulations}

The VLMC corresponding to the partition on example  (\label{exmodel}), have contexts:
\begin{eqnarray*}
T_1&=&\{0 \},\\
T_2&=&\{1\},\\
T_3&=&\{12 \},\\
T_4&=&\{102\},\\
T_5&=&\{202\},\\
T_6&=&\{22\},\\
T_7&=&\{002\}.
\end{eqnarray*}  

We simulated $1000$ samples of the process for each of the sample sizes $4000, 6000, 8000$ and $10000$. Using the tree as a basic good partition, for each sample we calculate the values $d_n(L_i,L_j)$ corresponding to the algorithm (\ref{al2}) and build the corresponding partitions. Table (\ref{tablaerr2}) show the results.

\begin{table}[htdp]

\caption{proportion of errors on the partition estimated for the model of example (\ref{exmodel})}\label{tablaerr2}
\begin{center}
\begin{tabular}{c|c}
Sample size & Proportion of errors\\
\hline
4000 & 0.614\\
6000 & 0.206\\
8000 & 0.047\\
10000 & 0.007\\
\end{tabular}
\end{center}
\label{default}
\end{table}%

Starting from the good partition corresponding to the context tree, the number of possible models is substantially reduced compared to those in the simulation on section (\ref{simulacion1}) and because of that, the error rates on this simulation are much better than before.

\section{Conclusions}\label{conclu} 

Our main motivation to define the minimal Markov models  is, in the first place, the concept of partitioning the state space in classes in which the states are equivalent, this allow us to model the  redundancy that appears in many processes  in the nature  as in genetics, linguistics, etc. Each class in the state space has a very specific, clear and practical meaning: any sequence of symbol in the same class has the same effect on the future distribution of the process. In other words, they activate the same random mechanism to choose the next symbol on the process. We can think of the resulting minimal partition as a  list of the relevant contexts for the process and their synonymous. 

In second place our motivation for developing this methodology is to demonstrate that  for a stationary, finite memory process it is theoretically  possible to find consistently a minimal Markov model to represent this process and that this can be accomplished in practice. The utilitarian implication of the fact that the model selection process can be started from a context tree partition, is that minimal Markov models can be easily fitted to stationary sources where the VLMC models already works.

It is clear that there are applications on which the natural partition to estimate is neither the minimal nor  a context tree partition. As long as the partition particular properties are well defined,  we can use theorem \ref{GeralREFTEO1} to estimate the minimal partition satisfying those properties. 

Our theorems are still valid if we change the constant term in the penalization of the BIC criterion for any positive (and finite) number. 
In the case of the VLMC model, the problem of finding a better constant has been addressed in diverse works as for example  \citet{Buhlmann1999} and \citet{galves2009}.

\section{Proofs} \label{demostra}

\begin{definition}\label{entropy} Let be $P$ and $Q$ probability distributions on $A.$ The relative entropy between $P$ and $Q$ is given by,
\begin{eqnarray}
D(P(\cdot) || Q(\cdot))= \sum_{a \in A} P(a) \ln \left( \frac{P(a)}{Q(a)} \right). \nonumber
\end{eqnarray} 
\end{definition}

\subsection{Proof of theorem  \ref{GeralREFTEO1}} 

\begin{eqnarray*}
BIC({\cal L}, x_1^n)
&=&\sum_{a \in A}\ln\left(\prod_{ L \in {\cal L}} \left(  \frac{r_n(L,a)}{r_n(L)} \right)^{N_n^{\cal L}(L,a)}\right)-\frac{(|{A}|-1)|{\cal L}|}{2}\ln(n),
\end{eqnarray*}
as consequence,
\begin{eqnarray}
BIC({\cal L}, x_1^n)-BIC({{\cal L}^{ij}}, x_1^n)&=& \sum_{ a \in A}\left\{
N^{{\cal L}}_n(L_i,a)\ln\left(\frac{r_n(L_i,a)}{r_n(L_i)}\right)\right. \nonumber \\ 
&&+ N^{{\cal L}}_n(L_j,a)\ln\left(\frac{r_n(L_j,a)}{r_n(L_j)}\right) \nonumber \\
&&-\left.  N^{{\cal L}^{ij}}_n(L_{ij},a)\ln\left(\frac{r_n(L_{ij},a)}{r_n(L_{ij})}\right)\right\}-\frac{(|{A}|-1)}{2}\ln(n). \label{dgeral}
\end{eqnarray}
We note that, the condition $\,\,I_{\{BIC({\cal L}^{ij}, x_1^n)>BIC({\cal L}, x_1^n)\}}=1$ is true if, and only if 
\begin{eqnarray} \sum_{ a \in A}\left\{r_n(L_i,a)\ln \left( \frac {r_n(L_i,a)}{r_n(L_i)}\right)
+r_n(L_j,a)\ln \left( \frac {r_n(L_j,a)}{r_n(L_j)}\right) \right. \nonumber \\
\left.
-r_n(L_{ij},a)\ln\left(\frac {r_n(L_{ij},a)}{r_n(L_{ij})}
\right) \right\}< \frac{(|{A}|-1)\ln(n)}{2n}.
\label{Ind1}
\end{eqnarray}
Because $r_n(L,a)$ and $r_n(L)$ are non-negative, using Jensen we have that,
\begin{eqnarray}
r_n(L_{i},a)\ln \left( \frac {r_n(L_{i},a)}{r_n(L_{i})}\right)
+r_n(L_{j},a)  \ln \left( \frac {r_n(L_{j},a)}{r_n(L_{j})} \right)   \geq \nonumber \\
 \left(r_n(L_{i},a)+ r_n(L_{j},a)  \right) \ln \left( \frac {r_n(L_{i},a)+ r_n(L_{j},a)}{r_n(L_{i})+r_n(L_{j})} \right) \nonumber
\end{eqnarray}
or equivalently,
\begin{equation} 
r_n(L_{i},a)\ln \left( \frac {r_n(L_{i},a)}{r_n(L_{i})} \right)
+r_n(L_{j},a)\ln \left( \frac {r_n(L_{j},a)}{r_n(L_{j})} \right) \geq r_n(L_{ij},a)\ln \left( \frac {r_n(L_{ij},a)}{r_n(L_{ij})} \right),
\label{Jensen}
\end{equation}
with equality if and only if  $\frac {r_n(L_{i},a)}{r_n(L_{i})}=\frac {r_n(L_{j},a)}{r_n(L_{j})},\,\, \forall a \in A.$\\
As consequence, equation (\ref{Jensen}) $\Rightarrow$ 
\begin{eqnarray}\label{Jensenconseq}
 \sum_{ a \in A}\left\{r_n(L_{i},a)\ln \left( \frac {r_n(L_{i},a)}{r_n(L_{i})}
 \right) +r_n(L_{j},a)\ln \left( \frac {r_n(L_{j},a)}{r_n(L_{j})}\right) \right. \nonumber \\
\left. -r_n(L_{ij},a)  \ln \left( \frac {r_n(L_{ij},a)}{r_n(L_{ij})} \right)
\right\} \geq 0, 
\end{eqnarray}
 with equality if and only if  $\frac {r_n(L_{i},a)}{r_n(L_{i})}=\frac {r_n(L_{j},a)}{r_n(L_{j})} \; \forall a \in A.$\\

Considering that $\frac{(|{A}|-1)\ln(n)}{2n} \to 0,$ as $n \to \infty$ and from the equation (\ref{Ind1}), we have that if\\
$\lim_{n \to \infty} I_{\{BIC({\cal L}^{ij}, x_1^n)>BIC({\cal L}, x_1^n)\}}=1,$ then 
 \begin{eqnarray*}
 \lim_{n \to \infty} \sum_{ a \in A}\left\{r_n(L_{i},a)\ln \left( \frac {r_n(L_{i},a)}{r_n(L_{i})} \right)
+r_n(L_{j},a)\ln \left( \frac {r_n(L_{j},a)}{r_n(L_{j})} \right) \right.\\
\left.-r_n(L_{ij},a)\ln \left( \frac {r_n(L_{ij},a)}{r_n(L_{ij})} \right)
\right\} \leq  0,
\end{eqnarray*}
from equation (\ref{Jensenconseq}) and taking the limit inside the sum we obtain
$$ \sum_{ a \in A}\left\{P(L_{i},a)\ln \left( \frac {P(L_{i},a)}{P(L_{i})}\right)
+P(L_{j},a)\ln \left( \frac {P(L_{j},a)}{P(L_{j})} \right)
-P(L_{ij},a)\ln \left(  \frac {P(L_{ij},a)}{P(L_{ij})} \right)
\right\} =  0,$$
using Jensen again, this means that
$\frac {P(L_{i},a)}{P(L_{i})}=\frac {P(L_{j},a)}{P(L_{j})} \; \forall a \in A,$
or equivalently,
$P(a|L_i)=P(a|L_j) \; \forall a \in A.$\\

For the other half of the proof, suppose  that 
$P(a|L_i)=P(a|L_j) \; \forall a \in A,$
as a consequence we have that 
\begin{eqnarray}
P(a|L_{ij})=P(a|L_i) \; \forall a \in A\label{paipb}	
\end{eqnarray}

\begin{eqnarray*}
BIC({\cal L}, x_1^n)-BIC({{\cal L}^{ij}}, x_1^n)&=&
	\ln\left(\prod_{ a \in A} \left(\frac{N^{{\cal L}}_n(L_{i},a)}{N^{{\cal L}}_n(L_{i})}\right)^{N^{{\cal L}}_n(L_{i},a)}\right)\\
	 &+& \ln\left(\prod_{ a \in A} \left( \frac{N^{{\cal L}}_n(L_{j},a)}{N^{{\cal L}}_n(L_{j})}\right)^{N^{{\cal L}}_n(L_{j},a)}\right)\\
	&-& \ln\left(\prod_{ a \in A} \left(\frac{N^{{\cal L}^{ij}}_n(L_{ij},a)}{N^{{\cal L}^{ij}}_n(L_{ij})}\right)^{N^{{\cal L}^{ij}}_n(L_{ij},a)}\right)-\frac{(|{A}|-1)}{2}\ln(n).
	\label{BIC1}
\end{eqnarray*}
Now, considering that  $\frac{N^{{\cal L}^{ij}}_n(L_{ij},a)}{N^{{\cal L}^{ij}}_n(L_{ij})}$ is the maximum likelihood estimator of $P(a|L_{ij})$,  
$$\prod_{ a \in A} \left(\frac{N^{{\cal L}^{ij}}_n(L_{ij},a)}{N^{{\cal L}^{ij}}_n(L_{ij})}\right)^{N^{{\cal L}^{ij}}_n(L_{ij},a)} \geq \prod_{ a \in A} P(a|L_{ij})^{N^{{\cal L}^{ij}}_n(L_{ij},a)}$$
$BIC({\cal L}, x_1^n)-BIC({{\cal L}^{ij}}, x_1^n)$ is bounded above by

\begin{eqnarray*}
&\ln&\left(\prod_{ a \in A} \left(\frac{N^{{\cal L}}_n(L_{i},a)}{N^{{\cal L}}_n(L_{i})}\right)^{N^{{\cal L}}_n(L_{i},a)}\right) + \ln\left(\prod_{ a \in A} \left(\frac{N^{{\cal L}}_n(L_{j},a)}{N^{{\cal L}}_n(L_{j})}\right)^{N^{{\cal L}}_n(L_{j},a)}\right)\\
&-&\ln\left(\prod_{ a \in A} P(a|L_{ij})^{N^{{\cal L}^{ij}}_n(L_{ij},a)}\right)-\frac{(|{A}|-1)}{2}\ln(n)\\
		&=&N^{{\cal L}}_n(L_{i})D\left(  \frac{N^{{\cal L}}_n(L_{i},.)}{N^{{\cal L}}_n(L_{i})} \Bigg| \Bigg| P(.|L_{i})\right)  +N^{{\cal L}}_n(L_{j})D\left(  \frac{N^{{\cal L}}_n(L_{j},.)}{N^{{\cal L}}_n(L_{j})} \Bigg| \Bigg| P(.|L_{j}) \right) - \frac{(|{A}|-1)}{2}\ln(n).
\end{eqnarray*}

Where $D(p||q)$ is the relative entropy, given by definition (\ref{entropy}). The first equality came from (\ref{paipb}) and (\ref{nlij}). Using proposition (\ref{lemaCsiszar06}), proposition  (\ref{lema2Csiszar02}), for any $\delta>0$ and $n$ large enough,
\begin{eqnarray}
D\left(  \frac{N^{{\cal L}}_n(L,.)}{N^{{\cal L}}_n(L)} \Bigg| \Bigg| P(.|L)\right) &\leq& \sum_{ a \in A} \frac{\left(  \frac{N^{{\cal L}}_n(L,a)}{N^{{\cal L}}_n(L)} - P(a|L)\right)^2}{P(a|L)}\\
 &\leq&  \sum_{ a \in A}\frac{\frac{\delta \ln(n)}{N^{{\cal L}}_n(L)}}{P(a|L)}.\label{RELENT1}
\end{eqnarray}

Then for any $\delta>0$ and $n$ large enough,  
\begin{eqnarray*} 
	BIC({\cal L}, x_1^n)-BIC({{\cal L}^{ij}}, x_1^n) &\leq& \frac{2\delta |A|}{p}\ln(n) - \frac{(|{A}|-1)}{2}\ln(n)\\
	 &=&\ln(n)\left( \frac{2\delta |A|}{p}-\frac{(|{A}|-1)}{2} \right)\\
\end{eqnarray*}
where $p=\min\{P(a|L): a\in A,L\in\{ L_i,L_j\} \}.$

 In particular, taking $\delta<\frac{p(|A|-1)}{4|A|}$, for $n$ large enough,
 $$BIC({\cal L}, x_1^n)-BIC({{\cal L}^{ij}}, x_1^n)<0.$$ 

\section*{Acknowledgements} We thank Antonio Galves, Nancy Garcia, Charlotte Galves and Florencia Leonardi for  their useful comments and discussions.


\begin{thebibliography}
{5}


\bibitem[Buhlmann P. and Wyner A. (1999)]{Buhlmann1999} \textsc{Buhlmann P.} and \textsc{Wyner A.} (1999). Variable length Markov chains. Ann. Statist. \textbf{27} 480--513. 
	\bibitem[Csisz\'{a}r, I. and Shields, P. C. (2000)]{Csiszar2000}  \textsc{Csisz\'{a}r, I. } and  \textsc{Shields, P. C.} (2000).The consistency of the BIC Markov order estimator. Ann. Statist. \textbf{28} 1601--1619. 
	\bibitem[Csisz\'{a}r, I. (2002)]{Csiszar2002} \textsc{Csisz\'{a}r, I. } (2002). Large-scale typicality of Markov sample paths and consistency of MDL order estimators. IEEE Trans. Inform. Theory \textbf{48} 1616--1628.
	\bibitem[Csisz\'{a}r, I. and Talata, Z. (2006)]{Csiszar2006} \textsc{Csisz\'{a}r, I.} and \textsc{Talata, Z.} (2006). Context tree estimation for not necessarily finite memory processes, via BIC and MDL. IEEE Trans. Inform. Theory \textbf{52} 1007--1016. 
	\bibitem[Galves, A., Galves, C., Garcia N. L. and Leonardi F. (2009)]{galves2009} \textsc{Galves et al.} (2009). Context tree selection and linguistic rhythm retrieval from written texts. arXiv:0902.3619. 
	\bibitem[Rissanen J. (1983)]{Rissanen1983} \textsc{Rissanen J.} (1983).  A universal data compression system,  IEEE Trans. Inform. 
Theory \textbf{29}(5) 656 -- 664. 
	\end{thebibliography}
\end{document}